\documentclass[12pt]{article}
\usepackage{amssymb,amsmath}
\usepackage{cases}
\usepackage{amsfonts}
\usepackage{graphicx}
\usepackage{cite,color,xcolor}
\usepackage[left=2.0cm,right=2.0cm,top=2.0cm,bottom=2.0cm]{geometry}
\usepackage[colorlinks,citecolor=blue,urlcolor=blue]{hyperref}
\usepackage[utf8]{inputenc}

\newtheorem{theorem}{Theorem}[section]
\newtheorem{corollary}{Corollary}[section]
\newtheorem{lemma}{Lemma}[section]
\newtheorem{proposition}{Proposition}[section]

\newtheorem{definition}{Definition}[section]
\newtheorem{remark}{Remark}[section]

\usepackage{txfonts}

\newcommand{\bal}{\begin{align}}
\newcommand{\bbal}{\begin{align*}}
\newcommand{\beq}{\begin{equation}}
\newcommand{\eeq}{\end{equation}}
\newcommand{\bca}{\begin{cases}}
\newcommand{\eca}{\end{cases}}

\newcommand{\pa}{\partial}
\newcommand{\fr}{\frac}
\newcommand{\na}{\nabla}
\newcommand{\De}{\Delta}
\newcommand{\de}{\delta}

\newcommand{\cd}{\cdot}
\newcommand{\ep}{\varepsilon}
\newcommand{\dd}{\mathrm{d}}
\newcommand{\B}{\dot{B}}

\newcommand{\R}{\mathbb{R}}

\newcommand{\D}{\mathrm{div}}

\newcommand{\ee}{\vec{e}}
\newcommand{\e}{A}
\begin{document}
\title{Ill-posedness for the stationary Navier-Stokes equations in critical Besov spaces}

\author{Jinlu Li$^{1}$, Yanghai Yu$^{2,}$\footnote{E-mail: lijinlu@gnnu.edu.cn; yuyanghai214@sina.com(Corresponding author); mathzwp2010@163.com} and Weipeng Zhu$^{3}$\\
\small $^1$ School of Mathematics and Computer Sciences, Gannan Normal University, Ganzhou 341000, China\\
\small $^2$ School of Mathematics and Statistics, Anhui Normal University, Wuhu 241002, China\\
\small $^3$ School of Mathematics and Big Data, Foshan University, Foshan, Guangdong 528000, China}

\date{\today}

\maketitle\noindent{\hrulefill}

{\bf Abstract:} This paper presents some progress toward an open question which proposed by Tsurumi (Arch. Ration. Mech. Anal. 234:2, 2019): whether or not the stationary Navier-Stokes equations in $\R^d$ is well-posed from $\dot{B}_{p, q}^{-2}$ to $\mathbb{P} \dot{B}_{p, q}^{0}$ with $p=d$ and $1 \leq q \leq 2$. In this paper, we prove that for the case $1\leq q<\frac d2$ with $d\geq4$ the stationary Navier-Stokes equations is ill-posed from $\dot{B}_{d, q}^{-2}(\R^d)$ to $\mathbb{P} \dot{B}_{d, q}^{0}(\R^d)$ by showing that a sequence of external forces is constructed to show discontinuity of the solution map at zero. Indeed in such case of $q$, there exists a sequence of external forces which converges to zero in $\dot{B}_{d, q}^{-2}$ and yields a sequence of solutions which does not converge to zero in $\dot{B}_{d, q}^{0}$. In particular, we also prove that the stationary Navier-Stokes equations is well-posed from $\dot{B}_{d, 2}^{-2}(\R^d)$ to $\mathbb{P} \dot{B}_{d, 2}^{0}(\R^d)$ with $d=3,4$. Based on these two cases, we demonstrate that the above open question for the dimension $d\geq4$ has been solved completely.

{\bf Keywords:} Stationary Navier-Stokes equations, Besov spaces, Ill-posedness

{\bf MSC (2010):} 35Q30; 35R25; 42B37
\vskip0mm\noindent{\hrulefill}

\section{Introduction}

In this paper, we consider the forced stationary Navier-Stokes equations describing the motion of incompressible fluid in the whole space $\R^d$, $d\geq3$
\begin{equation}\tag{SNS}
\begin{cases}
-\De u+u\cd\na u+\na \Pi=f,\\
\D u=0,
\end{cases}
\end{equation}
where $u=u(x)=\big(u^{1}(x), \ldots, u^{d}(x)\big)$ and $\Pi=\Pi(x)$ denote the unknown velocity vector and the unknown pressure at the point $x=\left(x_{1}, \ldots, x_{d}\right) \in \mathbb{R}^{d}$, respectively, while $f=f(x)=\big(f^{1}(x), \ldots, f^{d}(x)\big)$ denotes the given external force.

\subsection{Known Well/Ill-posedness (WP/IP) results}
There have been various studies on strong solutions $u$ to (SNS) for given data $f$
on the whole space $\R^d$. Leray \cite{le} and Ladyzhenskaya \cite{oa} showed the existence
of strong solutions to (SNS), and Heywood \cite{he} constructed solutions of (SNS) as a limit of
solutions of the non-stationary Navier-Stokes equations:
\begin{equation}\tag{\rm{NNS}}
\begin{cases}
\pa_tu-\De u+u\cd\na u+\na \Pi=0,\;& x\in \R^d,\ t>0, \\
\D u=0,\;& x\in \R^d,\ t>0,\\
u|_{t=0}=u_0,\;& x\in \R^d.
\end{cases}
\end{equation}
We should mention that Koch and Tataru \cite{Koch2001} obtained the global well-posedness of the 3D (NNS) for small initial data in the space $\mathrm{BMO}^{-1}=\dot{F}_{\infty, 2}^{-1}$. On the other hand, Bourgain and Pavlovi\'{c} \cite{Bou} showed the ill-posedness of (NNS) in $\dot{B}_{\infty,\infty}^{-1}$ (which includes $\mathrm{BMO}^{-1}$ ). Later on, the ill-posedness in $\dot{B}_{\infty, q}^{-1}$ with $1 \leq q<\infty$ was also showed by Yoneda in \cite{yon} $(2<q<\infty)$ and Wang in \cite{wang} $(1 \leq q \leq 2)$. These spaces play a crucial role since these are scaling invariant for the initial data $u_{0}$ in (NNS).

Chen \cite{chen} showed that for every small external force having a divergence-form $f = \D F$ with $F\in L^{d/2}(\R^d)$, there exists an
unique strong solution $u$ of (SNS) in $L^{d}(\R^d)$. Secchi \cite{sec} investigated existence and regularity of solutions to (SNS) in $L^d \cap L^p$ with $p > d$. As for the well-posedness of (SNS) in homogeneous Besov spaces, Kaneko-Kozono-Shimizu in \cite{Kan} showed the well-posed result as follows:

\begin{theorem}[see \cite{Kan}]\label{th1} Let $d \geq 3$. Suppose that $1 \leq p<d$ and $1 \leq q \leq \infty.$ Then $(\mathrm{SNS})$ is well-posed from $E=\dot{B}_{p, q}^{-3+\frac{d}{p}}$ to $S=\mathbb{P} \dot{B}_{p, q}^{-1+\frac{d}{p}}$.
\end{theorem}
Moreover, in the case $p=d$, (SNS) is also well-posed from $E=\dot{B}_{d, q}^{-2}$ to $S=\mathbb{P} L^{d}$ if $1 \leq q \leq 2$. These spaces $E$ and $S$ are scaling invariant for the external force $f$ and the velocity $u$ in (SNS) respectively. Precisely speaking, the corresponding scaling transform is $\{u, \pi, f\} \mapsto$ $\left\{u_{\lambda}, \Pi_{\lambda}, f_{\lambda}\right\}$ with $u_{\lambda}(x)=\lambda u(\lambda x), \Pi_{\lambda}(x)=\lambda^{2} \Pi(\lambda x), f_{\lambda}(x)=\lambda^{3} f(\lambda x)$, and we see that
$$
\left\|f_{\lambda}\right\|_{E}=\|f\|_{E}, \quad\left\|u_{\lambda}\right\|_{S}=\|u\|_{S}, \quad \forall \lambda>0.
$$
There are other previous results on the well-posedness in the case $p=d$. Bjorland et al. \cite{Bjorland} showed the well-posedness with more general space of external forces. In fact, they proved that there are constants $\varepsilon, \delta>0$ such that if $f \in \mathcal{S}^{\prime}$ satisfies $\left\|(-\Delta)^{-1} f\right\|_{L^{d, \infty}}<\varepsilon$, then there exists a unique solution $u \in B_{\mathbb{P} L^{d, \infty}}(\delta)$ to (SNS), which belongs to $L^{d}$ if and only if $\mathbb{P} f \in \dot{H}^{-2, d}$. Phan and Phuc \cite{Phan} showed the well-posedness in the largest critical space of external forces including $\dot{H}^{-2, d}$.

A nature question to ask is: Whether or not (SNS) is well-posed from $E=\dot{B}_{p, q}^{-3+\frac{d}{p}}$ to $S=\mathbb{P} \dot{B}_{p, q}^{-1+\frac{d}{p}}$ when $d\leq p \leq \infty$ and $1 \leq q \leq \infty$?

Recently, Tsuruni gave a partial answer to the above problem. More precisely, Tsuruni \cite{Tsu1} proved the ill-posedness of (SNS) in $\R^d$ (see \cite{Tsu3} for the Torus case $\mathbb{T}^d$), namely,

\begin{theorem}[see \cite{Tsu1}]\label{th2} Let $d<p \leq \infty, 1 \leq q \leq \infty$, and if $p=d, 2<q \leq \infty$, then (SNS) is ill-posed from $E=\dot{B}_{p, q}^{-3+\frac{d}{p}}$ to $S=$ $\mathbb{P} \dot{B}_{p, q}^{-1+\frac{d}{p}}$ in the sense that the solution map $f \in E \mapsto u \in S$ is, even if it exists, not continuous. More precisely, under such a condition, there exists a sequence $\left\{f_{N}\right\}_{N \in \mathbb{N}}$ of external forces with $f_{N} \rightarrow 0$ in $E$ such that there exists a unique solution $u_{N} \in \mathbb{P} L^{d}$ of (SNS) for each $f_{N}$, which never converges to zero in $S$ (actually, even in $\dot{B}_{\infty, \infty}^{-1}$ ).
\end{theorem}

Obviously, Tsuruni's result makes it clear that the well-posedness and ill-posedness can be divided between the case $(p, q) \in [1, d)
 \times [1,\infty]$ (Theorem \ref{th1}) and the case $(p, q) \in (d,\infty] \times [1,\infty]\cup \{d\} \times (2,\infty] $ (Theorem \ref{th2}), respectively. However, we remark that there is a gap between the global well-posedness in Theorem \ref{th1} and the ill-posedness in Theorem \ref{th2}. In other words, for the case $p=d$ and $1 \leq q \leq 2$, it is still unknown whether (SNS) in $\R^d$ is well-posed or ill-posed from $\dot{B}_{d, q}^{-2}$ to $\mathbb{P} \dot{B}_{d, q}^{0}$. In this paper, we are devoted to answer the question.
\subsection{New Results}
Now we return to the equation (SNS). Let us write the $i$-th component of $\mathbf{v}$ as $\mathbf{v}^{(i)}$. For the vector fields $\mathbf{v}$ and $\mathbf{u}$, we define the tensor product $\mathbf{v} \otimes \mathbf{u}$ as the $(i,j)$-th component $\mathbf{v}^{(i)} \mathbf{u}^{(j)}$ with $1 \leq i , j \leq d$. Thus, if $\mathbf{v}$ is divergence free (that is if $\nabla \cdot \mathbf{v}=0$ ) we have $\nabla \cdot(\mathbf{v} \otimes \mathbf{u})=(\mathbf{v} \cdot \nabla) \mathbf{u}$. Let us rewrite it to the generalized form so that we can apply successive approximation. First, we note that since $\D u=0$, there holds
$$
u \cdot \nabla u=\sum_{i=1}^{n} \partial_{ x_{i}}\left(u^{(i)} u\right)=\D(u \otimes u).
$$
We next introduce the projection $\mathbb{P}: L^{p}(\mathbb{R}^{d}) \rightarrow L_{\sigma}^{p}(\mathbb{R}^{d}) \equiv \overline{\left\{f \in \mathcal{C}^\infty_{0}(\mathbb{R}^{d}) ; \D f=0\right\}}^{\|\cdot\|_{L^{p}(\mathbb{R}^{d})}}$. In $\mathbb{R}^{d}$, $\mathbb{P}$ can be defined by $\mathbb{P}=\mathrm{Id}+\na(-\De)^{-1}\D$, or equivalently, $\mathbb{P}=(\mathbb{P}_{i j})_{1 \leqslant i, j \leqslant d}$ with $\mathbb{P}_{i j} \equiv \delta_{i j}+R_{i} R_{j}$.

Applying $\mathbb{P}$ to $(\mathrm{SNS})$, we obtain
$$
-\Delta u+\mathbb{P} \nabla \cdot(u \otimes u)=\mathbb{P}f,
$$
implied by $\mathbb{P} u=u$ and $\mathbb{P}(\nabla \Pi)=0$ since $\D u=0$.
Also, since $\mathbb{P}$ commutes with $-\Delta$, the solution $u$ of $(\mathrm{SNS})$ can be expressed as
\begin{align}
u &=(-\Delta)^{-1} \mathbb{P}(u \cdot \nabla u)+(-\Delta)^{-1}  \mathbb{P} f \nonumber\\
&=\mathbb{P}(-\Delta)^{-1} \D(u \otimes u)+\mathbb{P}(-\Delta)^{-1} f \nonumber\\
&=\mathcal{B}(u,u)+g,\tag{\rm{rSNS}}
\end{align}
here and in what follows, we shall denote the bilinear form
$$\mathcal{B}(u,v)\equiv\mathbb{P}(-\Delta)^{-1} \D(u \otimes v).$$

Our main results now read as follows:
\begin{theorem}\label{t3} Let $d=3,4$. $(\mathrm{SNS})$ is well-posed from $E=\dot{B}_{d, 2}^{-2}(\mathbb{R}^{d})$ to $S=\mathbb{P} \dot{B}_{d, 2}^{0}(\mathbb{R}^{d})$.
\end{theorem}

\begin{theorem}\label{th3}
Let $d\geq 4$ and $1\leq q<\frac{d}2$. (rSNS) is ill-posed from $\B^{0}_{d,q}(\R^d)$ to $\mathbb{P}\B^{0}_{d,q}(\R^d)$ in the following sense:
There exist $\left\{g_{n}\right\}_{n=1}^{\infty} \subset\dot{B}_{d, q}^{0}(\mathbb{R}^{d})$ such that a sequence of solutions $u_n\in \dot{B}_{d,q}^{0}(\mathbb{R}^{d})$ to (rSNS) which satisfies
$$
\lim _{n \rightarrow \infty}\left\|g_{n}\right\|_{\dot{B}_{d, q}^{0}(\mathbb{R}^{d})} = 0
$$
and for some positive constant $\ep_0$
$$
\lim _{n \rightarrow \infty}\left\|u_{n}\right\|_{{\B}_{d, q}^{0}(\mathbb{R}^{d})} \geq \ep_0.
$$
\end{theorem}
\begin{remark}
Theorem \ref{th3} demonstrates that if $d\geq 4$ and $1\leq q<\frac{d}2$, there exists a sequence of external forces which converges to zero in $\dot{B}_{d, q}^{0}$ and yields a sequence of solutions to (rSNS) which does not converge to zero in $\dot{B}_{d, q}^{0}$. In other words, the (SNS) is ill-posed from $\dot{B}_{d, q}^{-2}$ to $\mathbb{P} \dot{B}_{d, q}^{0}$ due to the discontinuity of the solution map at zero.
\end{remark}
\begin{remark}
We should mention that we have completely solved the open question which proposed by Tsurumi \cite{Tsu1} for the case $d\geq4$. This can be seen clearly from the Table below.
\begin{table}[http]
      \centering
      \begin{tabular}{|l|c|c|c|r|}\hline
        $d$ &$p$&$q$&$\dot{B}_{p, q}^{-3+d/p}\mapsto\mathbb{P} \dot{B}_{p, q}^{-1+d/p}$\\\hline
        $d\geq3$&$[1,d)$&$[1,\infty]$&WP, \;see Theorem \ref{th1}\\\hline
        $d\geq3$&$(d,\infty]$&$[1,\infty]$&IP, \;see Theorem \ref{th2}\\\hline
        $d\geq3$&$d$&$(2,\infty]$&IP, \;see Theorem \ref{th2}
        \\\hline
        $d=3,4$&$d$&$2$&WP, \;see Theorem \ref{t3}
          \\\hline
        $d=4$&$d$&$[1,2)$&WP, \;see Theorem \ref{th3}
        \\\hline
        $d\geq5$&$d$&$[1,2]$&IP, \;see Theorem \ref{th3}\\\hline
        $d=3$&$d$&$[1,2)$& Unknown \\\hline
        \end{tabular}
        \caption{Well/Ill-posedness}
        \end{table}
\end{remark}
\subsection{Main Ideas}

To illustrate our main idea, we introduce
\begin{equation}\label{g}
\begin{cases}
G\equiv\mathcal{B}(g,g), \\
U\equiv\mathcal{B}(u,u)-G,
\end{cases}
\end{equation}
then we have form (rSNS)
\begin{equation}\label{g1}
u=g+G+U,
\end{equation}
Obviously,
\bbal
U&=\mathcal{B}(U,g)+\mathcal{B}(g,U)+\mathcal{B}(U,U)
+\mathcal{B}(U,G)\\&\quad+\mathcal{B}(G,U)
+\mathcal{B}(g,G)+\mathcal{B}(G,g)+\mathcal{B}(G,G).
\end{align*}

\begin{lemma}\label{lyz} Let $d=3,4$. Then for $u, v \in \B^{0}_{d,2}(\mathbb{R}^{d})$, we have $\mathcal{B}(u,v) \in L^d(\mathbb{R}^{d})\cap \B^{0}_{d,2}(\mathbb{R}^{d})$ with the estimate
\bbal
\|\mathcal{B}(u,v)\|_{\B^{0}_{d,2}(\mathbb{R}^{d})}\leq  C\|u\|_{\B^{0}_{d,2}(\mathbb{R}^{d})}\|v\|_{\B^{0}_{d,2}(\mathbb{R}^{d})},
\end{align*}
where $C$ is a positive constant.
\end{lemma}
{\bf Proof.}\; Due to the embedding $L^{\frac d2}(\mathbb{R}^{d})\hookrightarrow\B^{0}_{\frac d2,2}(\mathbb{R}^{d})\hookrightarrow\B^{-1}_{d,2}(\mathbb{R}^{d})$ and $\B^{0}_{d,2}(\mathbb{R}^{d}\hookrightarrow L^d(\mathbb{R}^{d})$, we have
\bbal
\|\mathcal{B}(u,v)\|_{\B^{0}_{d,2}(\mathbb{R}^{d})}&\leq
 C\|u\otimes v\|_{L^{\frac d2}(\mathbb{R}^{d})}
\leq C\|u\|_{L^d(\mathbb{R}^{d})}\|v\|_{L^d(\mathbb{R}^{d})}
\leq C\|u\|_{\B^{0}_{d,2}(\mathbb{R}^{d})}\|v\|_{\B^{0}_{d,2}(\mathbb{R}^{d})}.
\end{align*}

\begin{remark} With Lemma \ref{lyz} at our disposal, we can prove Theorem \ref{t3} by successive approximation method. Since the procedure is standard (see \cite{Kan}), we shall not go into details.
\end{remark}
\begin{lemma}\label{l0} Let $d\geq 4$. Then for $u, v \in L^d(\mathbb{R}^{d})$, we have $\mathcal{B}(u,v) \in L^d(\mathbb{R}^{d})\cap \B^{0}_{d,\frac d2}(\mathbb{R}^{d})$ with the estimate
\bbal
\|\mathcal{B}(u,v)\|_{L^d(\mathbb{R}^{d})}+\|\mathcal{B}(u,v)\|_{\B^{0}_{d,\frac d2}(\mathbb{R}^{d})}\leq  C\|u\|_{L^d(\mathbb{R}^{d})}\|v\|_{L^d(\mathbb{R}^{d})},
\end{align*}
where $C$ is a positive constant.
\end{lemma}
{\bf Proof.}\; We have
\bbal
&\|\mathcal{B}(u,v)\|_{L^d(\mathbb{R}^{d})}\leq C\|u\otimes v\|_{H^{-1,d}(\mathbb{R}^{d})}\leq C\|u\otimes v\|_{L^{\frac d2}(\mathbb{R}^{d})}\leq C\|u\|_{L^d(\mathbb{R}^{d})}\|v\|_{L^d(\mathbb{R}^{d})},
\\&\|\mathcal{B}(u,v)\|_{\B^{0}_{d,\frac d2}(\mathbb{R}^{d})}\leq C\|u\otimes v\|_{\B^{-1}_{d,\frac d2}(\mathbb{R}^{d})}\leq C\|u\otimes v\|_{\B^{0}_{\frac d2,\frac d2}(\mathbb{R}^{d})}\leq C\|u\otimes v\|_{L^{\frac d2}(\mathbb{R}^{d})}\leq C\|u\|_{L^d(\mathbb{R}^{d})}\|v\|_{L^d(\mathbb{R}^{d})}.
\end{align*}
We should emphasize that, Lemma \ref{l0} tells us that if there exists small enough $\delta$ such that $\|g\|_{L^d(\mathbb{R}^{d})}<\delta$, then we have a unique solution $u\in L^d(\mathbb{R}^{d})$  to (rSNS) such that $\|u\|_{L^d(\mathbb{R}^{d})}<\delta$, which implies that the solution map $L^d(\mathbb{R}^{d}) \ni g  \mapsto u \in L^d(\mathbb{R}^{d})$ of (rSNS) is well-posed.
\begin{corollary} \label{c1}
Let $d\geq 4$ and $1\leq q<\frac{d}2$. Assume that $g\in L^d(\mathbb{R}^{d})\cap \B^{1}_{\frac d2,q}(\mathbb{R}^{d})$ and $\|g\|_{L^d(\mathbb{R}^{d})\cap \B^{0}_{d,q}(\mathbb{R}^{d})}<\de$ for small enough $\delta$, then there exists a solution $u\in L^d(\mathbb{R}^{d})\cap \B^{0}_{d,\frac d2}(\mathbb{R}^{d})$  to (rSNS). Moreover, $u\in \B^{0}_{d,q}(\mathbb{R}^{d})$.
\end{corollary} 
{\bf Proof.}\; Due to $\B^{1}_{\frac d2,q}(\mathbb{R}^{d}\hookrightarrow \B^{0}_{d,q}(\mathbb{R}^{d})\hookrightarrow\B^{0}_{d,\frac{d}{2}}(\mathbb{R}^{d})$, using Lemma \ref{l0} yields that $u\in L^d(\mathbb{R}^{d})\cap \B^{0}_{d,\frac d2}(\mathbb{R}^{d})$.
It is easy to deduce that
\bbal
\|\mathcal{B}(u,u)\|_{\B^{1}_{\frac d2,q}(\mathbb{R}^{d})}&\leq C\|u\otimes u\|_{\B^{0}_{\frac d2,q}(\mathbb{R}^{d})}
\leq C\|u\|_{\B^{1}_{\frac d2,q}(\mathbb{R}^{d})}\|u\|_{L^d(\mathbb{R}^{d})},
\end{align*}
from which and (rSNS), we obtain 
\bbal
\|u\|_{\B^{1}_{\frac d2,q}(\mathbb{R}^{d})}&\leq C\|u\|_{\B^{1}_{\frac d2,q}(\mathbb{R}^{d})} \|u\|_{L^d(\mathbb{R}^{d})}+\|g\|_{\B^{1}_{\frac d2,q}(\mathbb{R}^{d})}\leq C\delta\|u\|_{\B^{1}_{\frac d2,q}(\mathbb{R}^{d})} +\|g\|_{\B^{1}_{\frac d2,q}(\mathbb{R}^{d})},
\end{align*}
which implies $u\in \B^{1}_{\frac d2,q}(\mathbb{R}^{d})\hookrightarrow \B^{0}_{d,q}(\mathbb{R}^{d})$.

Furthermore, we have
\begin{corollary} \label{c2} Let $d\geq 4$. For small enough $\delta$ such that $\|g\|_{L^d(\mathbb{R}^{d})}<\delta$, then we have
\bbal
\|U\|_{L^d(\mathbb{R}^{d})}+\|U\|_{\B^0_{d,\frac{d}{2}}(\mathbb{R}^{d})}\leq  C\|g\|_{L^d(\mathbb{R}^{d})}^3,
\end{align*}
where $C$ is a positive constant.
\end{corollary} 
{\bf Proof.}\;
Using Lemma \ref{l0} yields
\bbal
\|U\|_{L^d(\mathbb{R}^{d})}+\|U\|_{\B^0_{d,\frac{d}{2}}(\mathbb{R}^{d})}&\leq C\left(\|U\|_{L^d(\mathbb{R}^{d})}\|g\|_{L^d(\mathbb{R}^{d})}
+\|U\|_{L^d(\mathbb{R}^{d})}^2+\|g\|_{L^d(\mathbb{R}^{d})}
\|G\|_{L^d(\mathbb{R}^{d})}+\|G\|_{L^d(\mathbb{R}^{d})}^2\right)\\
&\leq C\left(\|U\|_{L^d(\mathbb{R}^{d})}\|g\|_{L^d(\mathbb{R}^{d})}
+\|U\|_{L^d(\mathbb{R}^{d})}^2+\|g\|_{L^d(\mathbb{R}^{d})}^3
+\|g\|_{L^d(\mathbb{R}^{d})}^4\right),
\end{align*}
which enables us to complete the proof of Corollary \ref{c2}.

\begin{remark} From \eqref{g1} and Corollary \ref{c2}, we expect that the primarily affect which leads to the discontinuity to the solution of (rSNS) is that the worst term from nonlinear interactions $\mathcal{B}(g,g)$.
\end{remark}

\section{Preliminaries}
Firstly, let us recall that for all $f\in \mathcal{S}'$, the Fourier transform $\widehat{f}$, is defined by
$$
(\mathcal{F} f)(\xi)=\widehat{f}(\xi)=\int_{\R^d}e^{-ix\xi}f(x)\dd x \quad\text{for any}\; \xi\in\R^d.
$$
The inverse Fourier transform of any $g$ is given by
$$
(\mathcal{F}^{-1} g)(x)=\check{g}(x)=\frac{1}{(2 \pi)^d} \int_{\R^d} g(\xi) e^{i x \cdot \xi} \dd \xi.
$$

Next, we will recall some facts about the Littlewood-Paley (L-P) decomposition, the homogeneous Besov spaces and their some useful properties.
\begin{proposition}[L-P decomposition, See \cite{B}] Let $\mathcal{B}:=\{\xi\in\mathbb{R}^d:|\xi|\leq 4/3\}$ and $\mathcal{C}:=\{\xi\in\mathbb{R}^d: 3/4\leq|\xi|\leq 8/3\}.$
Choose a radial, non-negative, smooth function $\chi:\R^d\mapsto [0,1]$ such that it is supported in $\mathcal{B}$ and $\chi\equiv1$ for $|\xi|\leq3/4$. Setting $\varphi(\xi):=\chi(\xi/2)-\chi(\xi)$, then we deduce that $\varphi$ is supported in $\mathcal{C}$ and $\varphi(\xi)\equiv 1$ for $4/3\leq |\xi|\leq 3/2$. Moreover,
\begin{align*}
&\sum_{j\in\mathbb{Z}}\varphi(2^{-j}\xi)=1, \quad \forall \;  \xi\in \R^d\setminus\{0\},\\
&\frac{1}{2} \leq \sum_{j \in \mathbb{Z}} \varphi^{2}(2^{-j} \xi) \leq 1,\quad \forall \;  \xi\in \R^d\setminus\{0\}.
\end{align*}
\end{proposition}
For every $u\in \mathcal{S'}(\mathbb{R}^d)$, the homogeneous dyadic blocks ${\dot{\Delta}}_j$ is defined as follows
\bbal
&\dot{\Delta}_ju=\varphi(2^{-j}D)u=\mathcal{F}^{-1}\big(\varphi(2^{-j}\cdot)\mathcal{F}u\big)
=2^{dj}\int_{\R^d}\check{\varphi}\big(2^{j}(x-y)\big)u(y)\dd y,\quad \forall j\in \mathbb{Z},\\
&\widetilde{\dot{\Delta}}_ju=\widetilde{\varphi}(2^{-j}D)u=\sum_{|k-j|\leq1}{\dot{\Delta}}_{k}u.
\end{align*}
Moreover, the dyadic blocks ${\dot{\Delta}}_j$ satisfies the property of almost orthogonality:
 \begin{eqnarray*}
\dot{\Delta}_j\dot{\Delta}_ku\equiv0 \quad \text{if} \quad |j-k|\geq2.
\end{eqnarray*}
In the homogeneous case, the following Littlewood-Paley decomposition makes sense
$$
u=\sum_{j\in \mathbb{Z}}\dot{\Delta}_ju\quad \text{for any}\;u\in \mathcal{S}'_h(\R^d),
$$
where $\mathcal{S}'_h$ is given by
\begin{eqnarray*}
\mathcal{S}'_h:=\Big\{u \in \mathcal{S'}(\mathbb{R}^{d}):\; \lim_{j\rightarrow-\infty}\|\chi(2^{-j}D)u\|_{L^{\infty}}=0 \Big\}.
\end{eqnarray*}
We turn to the definition of the Besov Spaces and norms which will come into play in our paper.
\begin{definition}[\cite{B}]
Let $s\in\mathbb{R}$ and $(p,q)\in[1, \infty]^2$. The homogeneous Besov space $\dot{B}^s_{p,q}(\R^d)$ consists of all tempered distribution $f$ such that
\begin{align*}
\dot{B}_{p,q}^{s}=\Big\{f \in \mathcal{S}'_h(\mathbb{R}^{d}):\; \|f\|_{\dot{B}_{p,q}^{s}(\mathbb{R}^{d})}< \infty \Big\},
\end{align*}
where
\begin{numcases}{\|f\|_{\dot{B}^{s}_{p,q}(\R^d)}:=}
\left(\sum_{j\in\mathbb{Z}}2^{sjr}\|\dot{\Delta}_jf\|^r_{L^p(\R^d)}\right)^{1/q}, &if $1\leq q<\infty$,\nonumber\\
\sup_{j\in\mathbb{Z}}2^{sj}\|\dot{\Delta}_jf\|_{L^p(\R^d)}, &if $q=\infty$.\nonumber
\end{numcases}
\end{definition}
\begin{remark} We point out that the following properties will be used in the sequel.
\begin{itemize}
  \item For the homogeneous Besov spaces, we have the embedding properties as follows:
$$
\dot{B}_{p, q_{1}}^{s} \hookrightarrow \dot{B}_{p, q_{2}}^{s},\; s \in \mathbb{R}, 1 \leq p \leq \infty, 1 \leq q_{1} \leq q_{2} \leq \infty
$$
and
$$
\dot{B}_{p_{1}, q}^{s_{1}} \hookrightarrow \dot{B}_{p_{2}, q}^{s_{2}},\;-\infty<s_{2} \leq s_{1}<\infty, 1 \leq q \leq \infty, 1 \leq p_{1} \leq p_{2} \leq \infty
$$
with $s_{1}-d / p_{1}=s_{2}-d / p_{2}$.
\item $\dot{B}_{p, 2}^{0}$ is continuously included in $L^p$ and $L^p$ is continuously included in $\dot{B}_{p, p}^{0}$, namely,
$$
\dot{B}_{p, 2}^{0} \hookrightarrow L^p\hookrightarrow \dot{B}_{p, p}^{0} ,\;2 \leq p< \infty.
$$
\item $\dot{B}_{p, p}^{0}$ is continuously included in $L^p$ and $L^p$ is continuously included in $\dot{B}_{p, 2}^{0}$, namely,
$$
\dot{B}_{p, p}^{0} \hookrightarrow L^p\hookrightarrow \dot{B}_{p, 2}^{0} ,\;1< p\leq 2.
$$
  \item The Riesz potential $(-\Delta)^{\frac{\alpha}{2}} f \equiv \mathcal{F}^{-1}\left(|\xi|^{\alpha} \widehat{f}(\xi)\right)$ with $\alpha \in \mathbb{R}$ gives an isomorphism from $\dot{B}_{p, q}^{s+\alpha}$ onto $\dot{B}_{p, q}^{s}$ for any $s \in \mathbb{R}$ and $1 \leq p, q \leq \infty$, which implies that
$$
\left\|(-\Delta)^{\frac{\alpha}{2}} f\right\|_{\dot{B}_{p, q}^{s}} \approx\|f\|_{\dot{B}_{p, q}^{s+\alpha}} .
$$
\end{itemize}
\end{remark}

\section{Proof of Theorem}
\subsection{Construction of initial data}
Letting $n \gg 1,$ we write
$$n\in 16\mathbb{N}=\left\{16,32,48,\cdots\right\}\quad\text{and}\quad\mathbb{N}(n)=\left\{k\in 8\mathbb{N}: \frac{n}4 \leq k\leq \frac{n}2\right\},$$
 Let $0<\varepsilon \ll 1$ ($\varepsilon$ will be chosen below, see Remark \ref{re2}). We define the matric $A$ whose $(i,j)$-th component $(A)_{ij}$ with $1 \leq i , j \leq d$ is given by
\bbal
\left(A\right)_{ij}=\bca
\varepsilon, \quad  1\leq i=j\leq 2,\\
1, \quad  3\leq i=j\leq d,\\
0, \quad else,
\eca
\quad\text{and}\quad
\ee=\frac{\sqrt{2}}{2}(1,1,\underbrace{0,\cdots,0}_{d-2}).
\end{align*}
Define a scalar function $\widehat{\theta}\in \mathcal{C}^\infty_0(\mathbb{R})$ with values in $[0,1]$ which satisfies
\bbal
\widehat{\theta}(\xi)=
\bca
1, \quad \mathrm{if} \ |\xi|\leq \frac{1}{200 d},\\
0, \quad \mathrm{if} \ |\xi|\geq \frac{1}{100 d}.
\eca
\end{align*}
Let
$$\phi(x)=\theta(x_1)\theta(x_2)\cdots\theta(x_{d-1})\theta(x_d)
\sin\left(\frac{17}{24}x_d\right).$$
To construct a sequence of external forces $\{f_n\}^\infty_{n=1}$, motivated by Wang in \cite{wang} and Iwabuchi-Ogawa in \cite{Iwa}, we firstly need to introduce
\bal
b_n&\equiv(-\De)^{-1}a_n=n^{-\frac {1}{2q}}\sum\limits_{k\in \mathbb{N}(n)}2^{k}\phi\left(2^{k}A(x-2^{2n+k}\ee)\right)
\sin\left(\frac{17}{12}2^{n}\ee\cdot x\right),\label{b}\\
c_n&\equiv\mathcal{F}^{-1}\left(\frac{\xi_2-\xi_1}{\xi_2}\widehat{b}_n\right). \label{c}
\end{align}
Obviously, both $b_n$ and $c_n$ are real scalar functions. Also, it holds
$$(\pa_1-\pa_2)b_n=-\pa_2c_n.$$
In particular, we should emphasize the following important fact
\bal\label{z2}
\mathrm{supp} \ \widehat{b_n}(\xi)&\subset  \left\{\xi\in\R^d: \ \frac{33}{24}2^{n}\leq |\xi|\leq \frac{35}{24}2^{n}\right\}.
\end{align}
For the details of proof, see Lemma \ref{A1} in Appendix A.

We set
\begin{equation}\label{h}
\begin{cases}
g_n^{(1)}\equiv b_n ,\\
g_n^{(2)}\equiv c_n-b_n,\\
g_n^{(3)}\equiv
\cdots\equiv
g_n^{(d)}\equiv0.
\end{cases}
\end{equation}
It is clearly seen that $g_n\in \mathcal{S}(\R^d)$ and $\D g_n=0$, which obviously implies that $\mathbb{P}g_n=g_n$.

Hence, we can define the sequence of external forces $\{f_n\}^\infty_{n=1}$ by $\{g_n\}^\infty_{n=1}$, precisely,  $$g_n\equiv(-\De)^{-1}\mathbb{P}f_n.$$

Notice that for all $j\in \mathbb{Z}$
\bbal
\varphi(2^{-j}\xi)\equiv 1\quad \text{for}\quad \xi\in  \mathcal{C}_j\equiv\left\{\xi\in\R^d: \ \frac{4}{3}2^{j}\leq |\xi|\leq \frac{3}{2}2^{j}\right\},
\end{align*}
and
\bbal
\widehat{\dot{\Delta}_jb_n}=\varphi(2^{-j}\cdot)\widehat{b_n},
\end{align*}
which implies
\bbal
\widehat{\dot{\Delta}_jb_n}=0, \quad j\neq n.
\end{align*}
thus,
\begin{numcases}
{\dot{\Delta}_j(b_n)=}
b_n, &if $j=n$,\nonumber\\
0, &otherwise.\nonumber
\end{numcases}
Similarly, the above also holds for $c_n$ and $g_n$.
Based on the observation, we have the following two Lemmas involving $b_n$ and $c_n$.

We should remark that, here and in what follows, the positive constants $c$ and $C$ whose value may vary from line to line, may depend on $\ep$ and $\phi$ but not $n$. The positive constants $\widetilde{c}$ and $\widetilde{C}$ whose value may vary from line to line, may depend on $\phi$ but not $n$ and $\ep$.
\begin{lemma}\label{lem1}
Let $b_n$ be defined by \eqref{b}. Then there holds
\bbal
\|b_n\|_{\dot{B}^{0}_{d,1}}\leq Cn^{\frac1d-\frac{1}{2q}}.
\end{align*}
\end{lemma}
{\bf Proof.} Since $\phi$ is a Schwartz function, we have
\bal\label{s}
|\phi(x)|+\sum^d_{i=1}|\pa_{x_i}\phi(x)|\leq C(1+|x|)^{-M}, \qquad  M\geq 100d.
\end{align}
It is easy to show that
\bal\label{a0}
n^{\frac{d}{2q}}\left\|b_n\right\|^d_{L^d}&\leq  \int_{\R^d}\sum\limits_{\ell_1,\ell_2,\cdots,\ell_d\in \mathbb{N}(n)}\frac{2^{(\ell_1+\ell_2+\cdots+\ell_d)}}{(1+2^{\ell_1}|A(x
+2^{2n+\ell_1}\ee)|)^M\cdots (1+2^{\ell_d}|A(x+2^{2n+\ell_d})\ee|)^M}\dd x
\nonumber\\&\leq  \sum_{\ell\in \mathbb{N}(n)}\int_{\R^d}\frac{2^{d\ell}}{(1+2^{\ell}|A(x+2^{2n+\ell}\ee)|)^{dM}}\dd x\nonumber\\
&\quad+ \sum\limits_{(\ell_1,\ell_2,\cdots,\ell_d)\in \Lambda}\int_{\R^d}\frac{2^{(\ell_1+\ell_2+\cdots+\ell_d)}}{(1+2^{\ell_1}|A(x
+2^{2n+\ell_1}\ee)|)^M\cdots (1+2^{\ell_d}|A(x+2^{2n+\ell_d}\ee)|)^M}\dd x
\nonumber\\&\equiv I_1+I_2,
\end{align}
where the set $\Lambda$ is defined by
$$
\Lambda=\left\{\left(\ell_{1}, \ldots, \ell_{d}\right) \in \mathbb{N}^d(n)  \mid \exists 1 \leq k, \ell \leq d \text { s.t. } \ell_{k} \neq \ell_{\ell} \right\} .
$$
For the term $I_1$, by direct computations, one has
\bal\label{es-I1}
I_1=\sum_{\ell\in \mathbb{N}(n)}\int_{\R^d}\frac{1}{(1+|Ax|)^{dM}}\dd x=\frac{1}{\varepsilon^2}\sum_{\ell\in \mathbb{N}(n)}\int_{\R^d}\frac{1}{(1+|x|)^{dM}}\dd x\leq \frac{C}{\varepsilon^2}n.
\end{align}
For the term $I_2$, we assume that $\ell_1< \ell_2$ without loss of generality, then $\ell_2- \ell_1\geq 4$.
\bbal
&\quad\int_{\R^d}\frac{1}{(1+2^{\ell_1}|A(x
+2^{2n+\ell_1}\ee)|)^M (1+2^{\ell_2}|A(x+2^{2n+\ell_2}\ee)|)^M}\dd x\\
&=\left(\int_{\mathbf{A}_{\ell_1}}+\int_{\mathbf{A}_{\ell_1}^c}\right)\frac{1}{(1+2^{\ell_1}|A(x
+2^{2n+\ell_1}\ee)|)^M (1+2^{\ell_2}|A(x+2^{2n+\ell_2}\ee)|)^M}\dd x,
\end{align*}
where we defined the set $A_{\ell_1}$ by
$$
\mathbf{A}_{\ell_1}\equiv\left\{x:\left|A(x+2^{2 n+\ell_1} \ee) \right| \leq \varepsilon2^{2n}\right\}.
$$
Thus
\bal\label{a1}
&\int_{\mathbf{A}_{\ell_1}^c}\frac{1}{(1+2^{\ell_1}|A(x
+2^{2n+\ell_1}\ee)|)^M (1+2^{\ell_2}|A(x+2^{2n+\ell_2}\ee)|)^M}\dd x\nonumber\\
\leq&~ C(\varepsilon2^{\ell_{1}} 2^{2 n})^{-M}\int_{\mathbf{A}_{\ell_1}^c}\frac{1}{(1+2^{\ell_2}|A(x+2^{2n+\ell_2}\ee)|)^M}\dd x\nonumber\\
\leq&~ C(\varepsilon2^{\ell_{1}} 2^{2 n})^{-M}\varepsilon^{-2}2^{-d\ell_2}.
\end{align}
It is easy to deduce that for $x\in \mathbf{A}_{\ell_1}$
\begin{align*}
2^{\ell_{2}}\left|A(x+2^{2 n+\ell_2} \ee)\right| & \geq 2^{\ell_{2}}\left|A(2^{2 n+\ell_2}-2^{2 n+\ell_1}) \ee\right|- \varepsilon2^{\ell_{2}} 2^{2 n} \geq  \varepsilon2^{\ell_{2}} 2^{2 n}.
\end{align*}
Similarly,
\bal\label{a2}
&\int_{\mathbf{A}_{\ell_1}}\frac{1}{(1+2^{\ell_1}|A(x
+2^{2n+\ell_1}\ee)|)^M (1+2^{\ell_2}|A(x+2^{2n+\ell_2}\ee)|)^M}\dd x\nonumber\\
\leq&~ (\varepsilon2^{\ell_{2}} 2^{2 n})^{-M}\int_{\mathbf{A}_{\ell_1}}\frac{1}{(1+2^{\ell_1}|A(x
+2^{2n+\ell_1}\ee)|)^M }\dd x\nonumber\\
\leq&~ C(\varepsilon2^{\ell_{2}} 2^{2 n})^{-M}\varepsilon^{-2}2^{-d\ell_1}.
\end{align}
We infer from \eqref{a1} and \eqref{a2} that
\bal\label{es-I2}
&I_2\leq C\ep^{-M-2} 2^{-2 Mn}\sum\limits_{(\ell_1,\ell_2,\cdots,\ell_d)\in \Lambda}(2^{-M\ell_{1}}2^{-d\ell_2}+2^{-M\ell_{2}}2^{-d\ell_1})2^{(\ell_1+\ell_2+\cdots+\ell_d)}\leq C2^{-Mn}.
\end{align}
Inserting \eqref{es-I1} and \eqref{es-I2} into \eqref{a0}, we have for large enough n
\bbal
&\|b_n\|_{\dot{B}^{0}_{d,1}}\approx\|b_n\|_{L^d}\leq C n^{\fr{1}{d}-\fr1{2q}}.
\end{align*}
This completes the proof of Lemma \ref{lem1}.

\begin{lemma}\label{lem2}
Let $c_n$ be defined by \eqref{c}. Then there holds
\bbal
\|c_n\|_{\dot{B}^{0}_{d,1}}\leq C2^{-\frac n2}.
\end{align*}
\end{lemma}
{\bf Proof.}\; By Hausdorff-Young's inequality, we have
\bbal
\|c_n\|_{L^d}&\leq C\left\|\frac{\xi_1-\xi_2}{\xi_2}\widehat{b}_n(\xi)\right\|_{L^{\frac{d}{d-1}}}\\
&\leq C\sum\limits_{k\in \mathbb{N}(n)}2^{k}\left\|\frac{\xi_1-\xi_2}{\xi_2}\Phi_k^{\pm,\pm}(\xi)\right\|_{L^{\frac{d}{d-1}}},
\end{align*}
where $\Phi_k^{\pm,\pm}$ is given in  Appendix.

Noticing that the support condition of $\Phi_k^{\pm,\pm}$ (see \eqref{f1} in Appendix), which implies that
$|\xi_1-\xi_2|\leq 2^k$ and $|\xi_2|\approx2^n$,
thus we obtain
\bbal
\|c_n\|_{L^d}&\leq C2^{-n}\sum\limits_{k\in \mathbb{N}(n)}2^{2k}\left\|\Phi_k^{\pm,\pm}(\xi)\right\|_{L^{\frac{d}{d-1}}}\\
&\leq C2^{-n}\sum\limits_{k\in \mathbb{N}(n)}2^{k}\\&\leq C2^{-\frac n2},
\end{align*}
where we have used the simple fact
\bbal
\left\|\Phi_k^{\pm,\pm}(\xi)\right\|_{L^{\frac{d}{d-1}}}\leq C2^{-k}.
\end{align*}
Due to $\|c_n\|_{\dot{B}^{0}_{d,1}}\approx\|c_n\|_{L^d}$, we obtain the desired result and finish the proof of Lemma \ref{lem2}.

Combing Lemma \ref{lem1} and Lemma \ref{lem2} yields
\begin{proposition}\label{pro1}
Let $g_n$ be defined by \eqref{h}. Then
\bbal
&\|g_n\|_{\dot{B}^{0}_{d,1}}\leq C n^{\fr1d-\fr1{2q}}.
\end{align*}
\end{proposition}
\begin{remark}
From Proposition \ref{pro1}, we know that $g_n\in L^d(\mathbb{R}^{d})\cap \B^{1}_{\frac d2,q}(\mathbb{R}^{d})$ and $\|g\|_{L^d(\mathbb{R}^{d})\cap \B^{0}_{d,q}(\mathbb{R}^{d})}\leq C n^{\fr1d-\fr1{2q}}$. Then Corollary \ref{c1} tells us that the solution map $g_n\mapsto u_n\in \B^{0}_{d,q}(\mathbb{R}^{d})$.
\end{remark}
The following proposition is crucial for the proof of the discontinuity of solutions.
\begin{proposition}\label{pro2}
Let $g_n$ be defined by \eqref{h}. If $\ep$ is small enough and $n$ is large enough, then there exists $c>0$ independent of $n$ such that
\bbal
\|\mathcal{B}(g_n,g_n)\|_{\B^{0}_{d,q}(\mathbb{N}(n))}\geq c.
\end{align*}
\end{proposition}
{\bf Proof.} Recalling that the definition of $\mathcal{B}(g_n,g_n)$ and $\D g_n=0$, we have
$$\|\mathcal{B}(g_n,g_n)\|_{\B^{0}_{d,q}(\mathbb{N}(n))}\approx
\|\mathbb{P}(g_n\cd\na g_n)\|_{\B^{-2}_{d,q}(\mathbb{N}(n))}.$$
Noticing that $(\pa_1-\pa_2)b_n=-\pa_2c_n$, we have
\bbal
(g_n\cd\na g_n)^{(1)}&=\frac12(\pa_1-\pa_2)[b^2_n]+c_n\pa_2b_n\\
&=(\pa_1-\pa_2)[b^2_n]+\pa_2(b_nc_n),\\
(g_n\cd\na g_n)^{(2)}&=\frac12(-\pa_1+\pa_2)[b^2_n]+b_n\pa_1c_n-b_n\pa_2c_n-c_n\pa_2b_n
+c_n\pa_2c_n
\\&=(-\pa_1+\pa_2)[b^2_n]+(\pa_1-\pa_2)(b_nc_n)+\pa_2[c^2_n]-\pa_2(b_nc_n),\\
(g_n\cd\na g_n)^{(i)}&=0, \quad i=3,\cdots,d.
\end{align*}
Then, we can rewrite
\bal\label{gn}
g_n\cd\na g_n=E_n+F_n,
\end{align}
where
\bbal
E_n^{(1)}&=(\pa_1-\pa_2)[b^2_n],\quad  F_n^{(1)}=\pa_2(b_nc_n),\\
E_n^{(2)}&=-(\pa_1-\pa_2)[b^2_n], \quad F_n^{(2)}=(\pa_1-\pa_2)(b_nc_n)+\pa_2[c^2_n]-\pa_2(b_nc_n),\\
E_n^{(i)}&=F_n^{(i)}=0, \quad i=3,\cdots,d.
\end{align*}
Then, from \eqref{gn}, we have
\bbal
\left\|\mathbb{P}(g_n\cd\na g_n)\right\|_{\B^{-2}_{d,q}(\mathbb{N}(n))}&\geq
\left\|\mathbb{P}(E_n)\right\|_{\B^{-2}_{d,q}(\mathbb{N}(n))}-\left\|\mathbb{P}(F_n)\right\|_{\B^{-2}_{d,q}(\mathbb{N}(n))}\\
&\geq
\left\|\mathbb{P}(E_n)\right\|_{\B^{-2}_{d,q}(\mathbb{N}(n))}-C2^{-\frac n4}
\\&\geq
\left\|\big(\mathbb{P}(E_n)\big)^{(1)}\right\|_{\B^{-2}_{d,q}(\mathbb{N}(n))}-C2^{-\frac n4},
\end{align*}
where we have used
\bbal
\left\|\mathbb{P}(F_n)\right\|_{\B^{-2}_{d,q}(\mathbb{N}(n))}&\leq \left\|F_n^{(1)}\right\|_{\B^{-2}_{d,q}(\mathbb{N}(n))}+\left\|F_n^{(2)}\right\|_{\B^{-2}_{d,q}(\mathbb{N}(n))}\\
&\leq C \left(\left\|b_nc_n\right\|_{\B^{-1}_{d,q}(\mathbb{N}(n))}+\left\|c^2_n\right\|_{\B^{-1}_{d,q}(\mathbb{N}(n))}\right)\\
&\leq C \left(\left\|b_nc_n\right\|_{\B^{0}_{\fr d2,q}(\mathbb{N}(n))}+\left\|c^2_n\right\|_{\B^{0}_{\fr d2,q}(\mathbb{N}(n))}\right)\\
&\leq C n^{\fr1q-\fr2d}\left(\left\|b_nc_n\right\|_{\B^{0}_{\fr d2,\fr d2}(\mathbb{N}(n))}+\left\|c^2_n\right\|_{\B^{0}_{\fr d2,\fr d2}(\mathbb{N}(n))}\right)\\
&\leq C n^{\fr1q-\fr2d}\left(\|b_n\|_{L^d}\|c_n\|_{L^d}+\|c_n\|^2_{L^d}\right)\\
&\leq C2^{-\frac n4}.
\end{align*}
Notice that $\mathbb{P}=\mathrm{Id}+\na(-\De)^{-1}\D$, then
\bbal
\big(\mathbb{P}(E_n)\big)^{(1)}&=(\pa_1-\pa_2)[b^2_n]+\pa_1(-\De)^{-1}\D E_n\\
&=(\pa_1-\pa_2)[b^2_n]+\pa_1(\pa_1-\pa_2)^2(-\De)^{-1}[b^2_n],
\end{align*}
namely,
\bal\label{hy}
\left\|\big(\mathbb{P}(E_n)\big)^{(1)}\right\|_{\B^{-2}_{d,q}(\mathbb{N}(n))}= \|(\pa_1-\pa_2)[b^2_n]+\pa_1(\pa_1-\pa_2)^2(-\De)^{-1}[b^2_n]\|_{\B^{-2}_{d,q}(\mathbb{N}(n))}.
\end{align}
Due to \eqref{b}, direct computations gives that
\bbal
b^2_n&=n^{-\frac{1}{q}}(H_{1}+H_{2}),
\end{align*}
where
\bal
H_{1}&\equiv\sum\limits_{k\in \mathbb{N}(n)}2^{2k}\phi^2\big(2^{k}\e(x-2^{2n+k}\ee)\big)
\sin^2\left(\frac{17}{12}2^{n}\ee\cdot x\right),\label{h1}\\
H_2&\equiv\sum\limits_{k,j\in \mathbb{N}(n),k\neq j}2^{k+j}\phi\big(2^{k}\e(x-2^{2n+k}\ee)\big)
\phi\big(2^{j}\e(x-2^{2n+j}\ee)\big)\sin^2\left(\frac{17}{12}2^{n}\ee\cdot x\right).\label{h2}
\end{align}
We should emphasize that, the following cancelation holds
$$\|(\pa_1-\pa_2)H_{2}+\pa_1(\pa_1-\pa_2)^2(-\De)^{-1}H_{2}\|_{\B^{-2}_{d,q}(\mathbb{N}(n))}=0,$$
which immediately comes from the fact $\dot{\Delta}_{\ell}H_2=0$ for $ \ell\in\mathbb{N}(n)$ (for the proof see Lemma \ref{A2}).

Thus, \eqref{hy} reduces to
\bal\label{hy1}
\left\|\big(\mathbb{P}(E_n)\big)^{(1)}\right\|^q_{\B^{-2}_{d,q}(\mathbb{N}(n))}= n^{-1}\|(\pa_1-\pa_2)H_{1}+\pa_1(\pa_1-\pa_2)^2(-\De)^{-1}H_{1}\|^q_{\B^{-2}_{d,q}(\mathbb{N}(n))}.
\end{align}
For $\ell\in \mathbb{N}(n)$, we decompose the term $H_1$ as
\bbal
H_{1}&=
2^{2\ell}\phi^2\big(2^{\ell}\e(x-2^{2n+\ell}\ee)\big)\sin^2\left(\frac{17}{12}2^{n}\ee\cdot x\right)\\
&\quad+\sum\limits_{k\in \mathbb{N}(n),k\neq \ell}2^{2k}\phi^2\big(2^{k}\e(x-2^{2n+k}\ee)\big)\sin^2\left(\frac{17}{12}2^{n}\ee\cdot x\right),
\end{align*}
using the simple fact $\sin^2\alpha=(1-\cos2\alpha)/2$, then we have
\bbal
\dot{\Delta}_\ell H_1&=\fr12\dot{\Delta}_\ell\left(2^{2\ell}\phi^2\big(2^{\ell}\e(x-2^{2n+\ell}\ee)\big)\right)\\
&\quad+\frac12\dot{\Delta}_\ell\left(\sum\limits_{k\in \mathbb{N}(n),k\neq \ell}2^{2k}\phi^2\big(2^{k}\e(x-2^{2n+k}\ee)\big)\right)
\\&\equiv\dot{\Delta}_\ell H_{1,1}+\dot{\Delta}_\ell H_{1,2}.
\end{align*}
Let us introduce the set $\mathbf{B}_{\ell}$ defined by
$$
\mathbf{B}_{\ell}\equiv\left\{x:\big|A(x-2^{2n+\ell} \ee)\big|\leq 2^{-\ell}\right\}.
$$
Then, we have
\bbal
&\quad \|(\pa_1-\pa_2)H_{1}+\pa_1(\pa_1-\pa_2)^2(-\De)^{-1}H_{1}\|^q_{\B^{-2}_{d,q}(\mathbb{N}(n))}
\\&\geq \sum_{\ell\in \mathbb{N}(n) }2^{-2q\ell}\|\dot{\Delta}_\ell (\pa_1-\pa_2)H_{1}+\dot{\Delta}_\ell \pa_1(\pa_1-\pa_2)^2(-\De)^{-1}H_{1}\|^q_{L^d(\R^d)}
\\&\geq \sum_{\ell\in \mathbb{N}(n) }2^{-2q\ell}\|\dot{\Delta}_\ell (\pa_1-\pa_2)H_{1}+\dot{\Delta}_\ell \pa_1(\pa_1-\pa_2)^2(-\De)^{-1}H_{1}\|^q_{L^d(\mathbf{B}_{\ell})}
\\&\geq\sum_{\ell\in \mathbb{N}(n) }2^{-2q\ell}\|\dot{\Delta}_\ell (\pa_1-\pa_2)H_{1,1}\|^q_{L^d(\mathbf{B}_{\ell})}-\sum_{\ell\in \mathbb{N}(n) }2^{-2q\ell}\|\dot{\Delta}_\ell \pa_1(\pa_1-\pa_2)^2(-\De)^{-1}H_{1,1}\|^q_{L^d(\mathbf{B}_{\ell})}
\\&\quad-\sum_{\ell\in \mathbb{N}(n) }2^{-2q\ell}\|\dot{\Delta}_\ell (\pa_1-\pa_2)H_{1,2}\|^q_{L^d(\mathbf{B}_{\ell})}-\sum_{\ell\in \mathbb{N}(n) }2^{-2q\ell}\|\dot{\Delta}_\ell \pa_1(\pa_1-\pa_2)^2(-\De)^{-1}H_{1,2}\|^q_{L^d(\mathbf{B}_{\ell})}\\
&\equiv\mathbf{J_1}-\mathbf{J_2}-\mathbf{J_3}-\mathbf{J_4}.
\end{align*}
Now, we will show that the term $\mathbf{J_1}$ contributes the main part.

{\bf Estimation of $\mathbf{J_1}$.}\; Denote
\bbal
h(x)&\equiv\frac12(\pa_1-\pa_2)(\phi^2)(x)\\
&=\frac12(\pa_1-\pa_2)\left(\theta^2(x_1)\theta^2(x_2)\right)\theta^2(x_3)\cdots\theta^2(x_{d-1})\theta^2(x_d)
\sin^2\left(\frac{17}{24}x_d\right)\\
&=h_1(x)+h_2(x),
\end{align*}
where
\bbal
&h_1(x)\equiv\frac14(\pa_1-\pa_2)\left(\theta^2(x_1)\theta^2(x_2)\right)\theta^2(x_3)\cdots\theta^2(x_{d-1})\theta^2(x_d),\\
&h_2(x)\equiv\frac14(\pa_1-\pa_2)\left(\theta^2(x_1)\theta^2(x_2)\right)\theta^2(x_3)\cdots\theta^2(x_{d-1})\theta^2(x_d)
\cos\left(\frac{17}{12}x_d\right),
\end{align*}
by direct computations, we get for $\ell\in \mathbb{N}(n)$,
\bbal
\dot{\Delta}_\ell((\pa_1-\pa_2)H_{1,1})&=\varepsilon2^{3\ell}
h_2\big(2^{\ell}\e(x-2^{2n+\ell}\ee)\big).
\end{align*}
Then by change of variables, we have
\bbal
\|\dot{\Delta}_\ell((\pa_1-\pa_2)H_{1,1})\|_{L^{d}(\mathbf{B}_{\ell})}
&\geq \varepsilon2^{3\ell}\left\|h_2\big(2^{\ell}\e(x-2^{2n+\ell}\ee)\big)\right\|_{L^d(\mathbf{B}_{\ell})}\\
&\geq \varepsilon^{1-\frac2d}2^{2\ell}\left\|h_2(y)\right\|_{L^d(|y|\leq 1)}\\
&\equiv \widetilde{c}\varepsilon^{1-\frac2d}2^{2\ell},
\end{align*}
thus
\bal\label{j1}
\mathbf{J_1}=\sum\limits_{\ell\in \mathbb{N}(n)}2^{-2\ell q }\left\|{\Delta}_\ell\left((\pa_1-\pa_2)H_{1,1}\right)\right\|^q_{L^d(\mathbf{B}_{\ell})}\geq \widetilde{c}\varepsilon^{(1-\frac2d)q} n.
\end{align}
{\bf Estimation of $\mathbf{J_2}$.}\;
By direct computations, we get
\bbal
&\quad\|\dot{\Delta}_\ell(\pa_1(\pa_1-\pa_2)^2
(-\Delta)^{-1}H_{1,1})\|_{L^{d}(\mathbf{B}_{\ell})}\\
&\leq \|\dot{\Delta}_\ell(\pa_1(\pa_1-\pa_2)^2
(-\Delta)^{-1}H_{1,1})\|_{L^{d}(\R^d)}\\
&\leq \widetilde{C}2^{-2\ell}\|\dot{\Delta}_\ell(\pa_1(\pa_1-\pa_2)^2H_{1,1})\|_{L^{d}(\R^d)}\\
&\leq \widetilde{C}\ep^32^{3\ell}\|[(\pa_1(\pa_1-\pa_2)^2\phi^2]
\big(2^{\ell}\e(y-2^{2n+\ell}\ee)\big)\|_{L^{d}(\R^d)}
\\
&\leq \widetilde{C}\ep^{3-\frac2d}2^{2\ell}\|(\pa_1(\pa_1-\pa_2)^2\phi^2\|_{L^{d}(\R^d)}
\\&\leq \widetilde{C}\ep^{3-\frac2d}2^{2\ell},
\end{align*}
where we have used
\bbal
\dot{\Delta}_\ell(\pa_1(\pa_1-\pa_2)^2H_{1,1})&=\varepsilon^32^{5\ell}
[(\pa_1(\pa_1-\pa_2)^2\phi^2]\big(2^{\ell}\e(y-2^{2n+\ell}\ee)\big).
\end{align*}
Therefore
\bal\label{j2}
\mathbf{J_2}=\sum\limits_{\ell\in \mathbb{N}(n)}2^{-2\ell q }\left\|\dot{\Delta}_\ell\left(\pa_1(\pa_1-\pa_2)^2
(-\Delta)^{-1}H_{1,1}\right)\right\|^q_{L^d(\mathbf{B}_{\ell})}\leq \widetilde{C}\varepsilon^{(3-\frac2d)q} n.
\end{align}
\begin{remark}\label{re2}
It is worthwhile pointing out that both the constants $\widetilde{c}$ and $\widetilde{C}$ do not depend on the parameter $\ep$. Combining \eqref{j1} and \eqref{j2}, we need to take $\ep$ sufficiently small such that $\mathbf{J_2}$ can be absorbed by $\mathbf{J_1}$. This is why in the present paper we set the parameter $\ep$ in the matric $\e$.
\end{remark}
Next, we deal with the last  two terms which are much less than the first one.
By making fully use of the information of decreasing rapidly functions $\check{\varphi}$ and $\phi$ when $k\neq \ell$, we can deal with the third term.

{\bf Estimation of $\mathbf{J_3}$.}\;
Noting the fact that for $100d<N\in \mathbb{Z}^+$
$$|\check{\varphi}(x)|\leq C(1+|x|)^{-N},$$
then we have
\bal\label{j-3}
&\quad\left\|\dot{\Delta}_\ell\left((\pa_1-\pa_2)H_{1,2}\right)\right\|_{L^d(\mathbf{B}_{\ell})}\nonumber\\
&\leq\sum\limits_{k\in \mathbb{N}(n),k\neq \ell}2^{3k}2^{d\ell}\left\|\int_{\R^d}\check{\varphi}(2^{\ell}(x-y))
[(\pa_1-\pa_2)\phi^2]\big(2^{k}\e(y-2^{2n+k}\ee)\big)
\dd y\right\|_{L^d(\mathbf{B}_{\ell})}\nonumber\\
&\leq  \sum\limits_{k\in \mathbb{N}(n),k\neq \ell}2^{3k}2^{d\ell}\left\|\int_{\R^d}\Big(1+2^{\ell}|x-y|\Big)^{-N}\Big(1+2^{k}|A(y-2^{2n+k}\ee)|\Big)^{-2N}\dd y\right\|_{L^d(\mathbf{B}_{\ell})}.
\end{align}
Dividing the integral region in terms of $y$ into the following two parts to estimate:
\begin{align*}
\R^d&=\left\{y :\; |A(y-2^{\ell+2 n} \ee )| \leq \varepsilon2^{2 n}\right\}\cup\left\{y :\; |A(y-2^{\ell+2 n} \ee)|\geq \varepsilon2^{2 n}\right\}\\
&\equiv \mathbf{A}_{1} \cup \mathbf{A}_{2},
\end{align*}
For $x \in \mathbf{B}_{\ell}$ and $y \in \mathbf{A}_{1}$, we conclude that
$$
\begin{aligned}
\left|A(y-2^{k+2 n} \ee)\right| &=\left|A(y-2^{\ell+2 n} \ee)+A(2^{\ell+2 n} \ee-2^{k+2 n} \ee)\right| \\
& \geq\left|(2^{\ell+2 n} -2^{k+2 n} )A\ee\right|-\left|A(y-2^{\ell+2 n} \ee)\right| \\
&\geq \varepsilon2^{2 n}.
\end{aligned}
$$
For $x \in \mathbf{B}_{\ell}$ and $y \in \mathbf{A}_{2}$, it is easy to check that
$$
\begin{aligned}
|x-y| &\geq\left|A(y-2^{\ell+2 n} \ee)\right|-\left|A(x-2^{\ell+2 n} \ee)\right| \geq\varepsilon2^{2 n}- 2^{-\ell} \geq \ep2^{2 n-1}.
\end{aligned}
$$
Then, we have
\bbal
&\left\|\int_{\R^d}\Big(1+2^{\ell}|x-y|\Big)^{-N}\Big(1+2^{k}|A(y-2^{2n+k}\ee)|\Big)^{-2N}\dd y\right\|_{L^d(\mathbf{B}_{\ell})}\\
\leq&~  C2^{-2(k+2n)N}\left\|\int_{\mathbf{A}_1}\Big(1+2^{\ell}|x-y|\Big)^{-N}\dd y\right\|_{L^d(\mathbf{B}_{\ell})}\\
&\quad+C2^{-(\ell+2 n)N}\left\|\int_{\mathbf{A}_2}\Big(1+2^{k}|A(y-2^{2n+k}\ee)|\Big)^{-2N}\dd y\right\|_{L^d(\mathbf{B}_{\ell})}\\
\leq&~  C\left(2^{-d\ell}2^{-2(k+2n)N}+2^{-(\ell+2 n)N}2^{-dk}\right)2^{-\frac{\ell}{d}}.
\end{align*}
Plugging the above into \eqref{j-3} yields
\bbal
\left\|{\dot{\Delta}}_\ell\left((\pa_1-\pa_2)H_{1,2}\right)
\right\|_{L^d(\mathbf{B}_{\ell})}\leq C\sum\limits_{k\in \mathbb{N}(n),k\neq \ell}2^{3k}2^{d\ell}\left(2^{-d\ell}2^{-2(k+2n)N}+2^{-(\ell+2 n)N}2^{-dk}\right)2^{-\frac{\ell}{d}}\leq C2^{-n}.
\end{align*}
Therefore, we have
\bbal
\mathbf{J_3}=\sum\limits_{\ell\in \mathbb{N}(n)}2^{-2\ell q }\left\|{\dot{\Delta}}_\ell\left((\pa_1-\pa_2)H_{1,2}\right)\right\|^q_{L^d(\mathbf{B}_{\ell})}\leq C n2^{-qn}.
\end{align*}
{\bf Estimation of $\mathbf{J_4}$.}\;
Notice that
\bbal
\phi^2(x)&=\prod_{i=1}^d\theta^2(x_i)+
\frac12\prod_{i=1}^d\theta^2(x_i)\cos\left(\frac{17}{12}x_d\right)\\
&\equiv\Theta_1(x)+\Theta_2(x),
\end{align*}
then we have for $k\neq \ell$,
\bbal
 \dot{\Delta}_{\ell}\left(\Theta_2\big(2^{k}\e(x-2^{2n+k}\ee)
\big)\right)=0,
\end{align*}
which implies
\bbal
 \dot{\Delta}_{\ell}\Big(\pa_1(\pa_1-\pa_2)^2(-\Delta)^{-1}[\phi^2\big(2^{k}\e(x-2^{2n+k}\ee)
\big)]\Big)&= \dot{\Delta}_{\ell}\Big(\pa_1(\pa_1-\pa_2)^2(-\Delta)^{-1}[\Theta_1\big(2^{k}\e(x-2^{2n+k}\ee)
\big)]\Big)
\\&= \dot{\Delta}_{\ell}\widetilde{\dot{\Delta}}_{\ell}\Psi(x)
\end{align*}
where we denote
$$\Psi(x)\equiv U\big(2^{k}\e(x-2^{2n+k}\ee)\big).$$
and
$$U(x)=\ep^3 2^k \pa_1(\pa_1-\pa_2)^2
(-\ep^2\pa^2_1-\ep^2\pa^2_2-\pa^2_3-\cdots-\pa^2_d)^{-1}[\Theta_1](x),$$
we can show that
\bbal
\widetilde{\dot{\Delta}}_{\ell}\Psi(x)&=2^{d\ell}\int_{\R^d}
\tilde{h}(2^\ell(x-y))U\big(2^{k}\e(y-2^{2n+k}\ee)\big)\dd y
\\&=2^{d\ell}\int_{\R^d}
\tilde{h}(2^\ell(x-z-2^{2n+k}\ee))U\big(2^{k}\e z\big)\dd z
\\&=\ep^{-2}2^{d(\ell-k)}\int_{\R^d}
\tilde{h}(2^\ell(x-2^{-k}A^{-1}y-2^{2n+k}\ee))U\big(y\big)\dd y
\\&=\ep^{-2}2^{d(\ell-k)}\int_{\R^d}
\tilde{h}(2^{\ell-k}A^{-1}(2^kA[x-2^{2n+k}\ee]-y))U\big(y\big)\dd y
\\&=(P_\ell*U)\big(2^{k}\e(x-2^{2n+k}\ee)\big),
\end{align*}
where
$${P_\ell*U}=\ep^{-2}2^{d(\ell-k)}\int_{\R^d}
\tilde{h}(2^{\ell-k}A^{-1}(x-y))U\big(y\big)\dd y.$$
Notice that
\bbal
\mathcal{F}[\tilde{h}(2^{\ell-k}A^{-1}\cdot)]=\ep^22^{-d(\ell-k)}\tilde{\varphi}(2^{k-\ell}A\xi),
\end{align*}
we have
\bbal
(P_\ell*U)(x)=(2\pi)^{-d}\ep^32^{k}\int_{\R^d}e^{ix\cdot\xi}\cdot
\frac{-i\xi_1(\xi_1-\xi_2)^2}{\ep^2\xi_1^2+\ep^2\xi_2^2+\xi_3^2+\cdots+\xi^2_d}
\widehat{\Theta_1}(\xi)\tilde{\varphi}(2^{k-\ell}A\xi)\dd \xi.
\end{align*}
Since $\tilde{\varphi}(x)=\varphi(\frac12x)+\varphi(x)+\varphi(2x)$, then we have $\tilde{\varphi}(x)$ is supported in $\mathcal{C}:=\{\xi\in\mathbb{R}^d: 3/8\leq|\xi|\leq 16/3\}.$ Therefor, if $|\xi|\leq \frac382^{\ell-k}$, we can deduce that $\tilde{\varphi}(2^{k-\ell}A\xi)=0$.
Then, using
\bbal
(1+|x|^2)^N(P_\ell*U)(x)&=(2\pi)^{-d}\ep^32^{k}\int_{\R^d}(1-\De_\xi)^N(e^{ix\cdot \xi})\frac{-i\xi_1(\xi_1-\xi_2)^2}{\ep^2\xi_1^2+\ep^2\xi_2^2+\xi_3^2+\cdots+\xi^2_d}
\widehat{\Theta_1}(\xi)\tilde{\varphi}(2^{k-\ell}A\xi)\dd \xi
\\&=(2\pi)^{-d}\ep^32^{k}\int_{\R^d}e^{ix\cdot \xi}(1-\De_\xi)^N\left(\frac{-i\xi_1(\xi_1-\xi_2)^2}{\ep^2\xi_1^2+\ep^2\xi_2^2+\xi_3^2+\cdots+\xi^2_d}
\widehat{\Theta_1}(\xi)\tilde{\varphi}(2^{k-\ell}A\xi)\right)\dd \xi,
\end{align*}
we have
\bbal
|(1+|x|^2)^N(P_\ell*U)(x)|\leq C2^{k}2^{(k-\ell)2N},
\end{align*}
which implies
\bbal
\Big|\widetilde{\dot{\Delta}}_{\ell}\Psi(x)\Big|
\leq C2^{k}2^{(k-\ell)2N}\Big(1+2^{k}|A(x-2^{2n+k}\ee)|\Big)^{-2N}.
\end{align*}
Then, we have
\bbal
&\quad \left\|{\dot{\Delta}}_\ell\left(\pa_1(\pa_1-\pa_2)^2(-\Delta)^{-1}H_{1,2}\right)
\right\|_{L^d(\mathbf{B}_{\ell})}
\\&\leq
\sum\limits_{k\in \mathbb{N}(n),k\neq \ell}2^{2k}2^{d\ell}\left\|\int_{\R^d}\check{\varphi}(2^{\ell}(x-y))
\Big(\widetilde{\dot{\Delta}}_{\ell}(\Psi(y))\Big)
\dd y\right\|_{L^d(\mathbf{B}_{\ell})}\\
&\leq  \sum\limits_{k\in \mathbb{N}(n),k\neq \ell}2^{3k}2^{d\ell}2^{(k-\ell)2N}\left\|\int_{\R^d}\Big(1+2^{\ell}|x-y|\Big)^{-N}\Big(1+2^{k}|\e(y-2^{2n+k}\ee)
\big|\Big)^{-2N}\dd y\right\|_{L^d(\mathbf{B}_{\ell})}\\
&\leq  C\sum\limits_{k\in \mathbb{N}(n),k\neq \ell}2^{3k}2^{d\ell}2^{(k-\ell)2N}\left(2^{-d\ell}2^{-2(k+2n)N}+2^{-(\ell+2 n)N}2^{-dk}\right)2^{-\frac{\ell}{d}}\\
&\leq C2^{-n},
\end{align*}
this gives that
\bbal
\mathbf{J_4}\leq Cn2^{-qn}.
\end{align*}
Combining the above estimates yields
\bbal
\left\|\big(\mathbb{P}(E_n)\big)^{(1)}\right\|^q_{\B^{-2}_{d,q}(\mathbb{N}(n))}&
\geq(\widetilde{c}-\widetilde{C}\ep^2)\varepsilon^{(1-\frac2d)q}-C2^{-qn},
\end{align*}
which implies
\bal\label{hy2}
\|\mathcal{B}(g_n,g_n)\|_{\B^{0}_{d,q}(\mathbb{N}(n))}\geq c.
\end{align}
Thus, we complete the proof of Proposition \ref{pro2}.

Combining Propositions \ref{pro1} and \ref{pro2}, we shall prove Theorem \ref{th3}.

{\bf Proof of Theorem \ref{th3}.}\; Letting $u_n=g_n+\mathcal{B}(g_n,g_n)+U_n$, then we deduce from Proposition \ref{pro1} that
\bbal
\|g_n\|_{\B^{0}_{d,q}(\mathbb{N}(n))}\leq C\|g_n\|_{\B^{0}_{d,1}}\leq  Cn^{\frac{1}{d}-\frac{1}{2q}}
\end{align*}
and Corollary \ref{c2} that
\bbal
\|U_n\|_{\B^{0}_{d,q}(\mathbb{N}(n))}\leq n^{\frac1q-\frac{2}{d}} \|U_n\|_{\B^{0}_{d,\frac d2}(\mathbb{N}(n))}\leq Cn^{\frac1q-\frac{2}{d}}\|g_n\|^3_{L^d}\leq Cn^{\frac1q-\frac{2}{d}}n^{\frac3d-\frac{3}{2q}}\leq Cn^{\frac1d-\frac{1}{2q}}.
\end{align*}
Therefore, we obtain for large $n$ enough
\bbal
\|u_n\|_{\B^{0}_{d,q}}&\geq \|u_n\|_{\B^{0}_{d,q}(\mathbb{N}(n))}
\\&\geq \|\mathcal{B}(g_n,g_n)\|_{\B^{0}_{d,q}(\mathbb{N}(n))}
-\|g_n\|_{\B^{0}_{d,q}(\mathbb{N}(n))}
-\|U_n\|_{\B^{0}_{d,q}(\mathbb{N}(n))}
\\&\geq c-Cn^{\frac1d-\frac 1{2q}}\geq \frac{c}{2}.
\end{align*}

Thus, we have completed the proof of Theorem \ref{th3}.

\section{Appendix}\label{sec4}
For the sake of convenience, here we present more details in the computations.
\begin{lemma}\label{A1}
Let $b_n$ be defined by \eqref{b}. Then there holds
\bbal
\mathrm{supp} \ \widehat{b_n}(\xi)\subset  \left\{\xi\in \R^d: \ \frac{33}{24}2^{n}\leq |\xi|\leq \frac{35}{24}2^{n}\right\}.
\end{align*}
\end{lemma}
{\bf Proof.}\; For the sake of simplicity, we denote
$$\Phi_k:=\phi\left(2^{k}A(x-2^{2n+k}\ee)\right)
\sin\left(\lambda_n\ee\cdot x\right)\quad\text{with}\quad \lambda_n:=\frac{17}{12}2^{n}.$$
Using the fact $\sin\alpha=\frac{i}{2} (e^{-i \alpha}-e^{i \alpha})$ from Euler's formula, we deduce easily that
\begin{align*}
\mathcal{F}\left(\Phi_k\right)
&=\int_{\R^d} e^{-i x\cdot \xi} \phi\left(2^{k}\e(x-2^{2n+k}\ee)\right)
\sin\left(\lambda_n\ee\cdot x\right) \dd x \\
&=\frac{i}{2} \int_{\R^d}\left(e^{-i x\cdot\left(\xi+\lambda_n\ee\right)}-e^{-i x\cdot\left(\xi-\lambda_n\ee\right)}\right) \phi\left(2^{k}\e(x-2^{2n+k}\ee)\right) \dd x \\
&=\frac{i}{2} \left(\Phi_k^1-\Phi_k^2\right),
\end{align*}
where
\begin{align*}
&\Phi_k^1:=\int_{\R^d}e^{-i (x+2^{2n+k}\ee)\cdot\left(\xi+\lambda_n\ee\right)} \phi\left(2^{k}\e x\right) \dd x,\\
&\Phi_k^2:=\int_{\R^d}e^{-i (x+2^{2n+k}\ee)\cdot\left(\xi-\lambda_n\ee\right)} \phi\left(2^{k}\e x\right) \dd x.
\end{align*}
Let $\widetilde{\lambda}_n:=\frac{\sqrt{2}}{2}\lambda_n$. Then by change of variables, we have
\begin{align*}
\Phi_k^1
=&~\frac{1}{\varepsilon^{2}2^{dk}}e^{-i 2^{2n+k}\ee\cdot\left(\xi+\widetilde{\lambda}_n\ee\right)}\int_{\R}e^{-i x_1\cdot\varepsilon^{-1}2^{-k}\left(\xi_1+\widetilde{\lambda}_n\right)} \theta(x_1) \dd x_1\int_{\R}e^{-i x_2\cdot\varepsilon^{-1}2^{-k}\left(\xi_2+\widetilde{\lambda}_n\right)} \theta(x_2) \dd x_2\\
&\qquad\times\int_{\R}e^{-i x_3\cdot2^{-k}\xi_3} \theta(x_3) \dd x_3\cdots\int_{\R}e^{-i x_{d-1}\cdot2^{-k}\xi_{d-1}} \theta(x_{d-1}) \dd x_{d-1}\int_{\R}e^{-i x_d\cdot2^{-k}\xi_d} \theta(x_d)\sin\left(\fr{17}{24}x_d\right) \dd x_d\\
=&~\frac{i}{\varepsilon^{2}2^{dk+1}}e^{-i 2^{2n+k}\ee\cdot\left(\xi+\widetilde{\lambda}_n\ee\right)}\widehat{\theta}\left(\frac{\xi_1+\widetilde{\lambda}_n}{\varepsilon2^{k}}\right)
\widehat{\theta}\left(\frac{\xi_2+\widetilde{\lambda}_n}{\varepsilon2^{k}}\right)\widehat{\theta}\left(\frac{\xi_3}{2^{k}}\right)\cdots
\widehat{\theta}\left(\frac{\xi_{d-1}}{2^{k}}\right)\left[\widehat{\theta}\left(\frac{\xi_d
}{2^k}+\frac{17}{24}\right)-\widehat{\theta}\left(\frac{\xi_d
}{2^k}-\frac{17}{24}\right)\right]\\
=&~\frac{i}{\varepsilon^{2}2^{dk+1}}e^{-i 2^{2n+k}\ee\cdot\left(\xi+\widetilde{\lambda}_n\ee\right)}\widehat{\theta}\left(\frac{\xi_1+\widetilde{\lambda}_n}{\varepsilon2^{k}}\right)
\widehat{\theta}\left(\frac{\xi_2+\widetilde{\lambda}_n}{\varepsilon2^{k}}\right)\widehat{\theta}\left(\frac{\xi_3}{2^{k}}\right)\cdots
\widehat{\theta}\left(\frac{\xi_{d-1}}{2^{k}}\right)\widehat{\theta}\left(\frac{\xi_d
}{2^k}+\frac{17}{24}\right)\\
-&~\frac{i}{\varepsilon^{2}2^{dk+1}}e^{-i 2^{2n+k}\ee\cdot\left(\xi+\widetilde{\lambda}_n\ee\right)}\widehat{\theta}\left(\frac{\xi_1+\widetilde{\lambda}_n}{\varepsilon2^{k}}\right)
\widehat{\theta}\left(\frac{\xi_2+\widetilde{\lambda}_n}{\varepsilon2^{k}}\right)\widehat{\theta}\left(\frac{\xi_3}{2^{k}}\right)\cdots
\widehat{\theta}\left(\frac{\xi_{d-1}}{2^{k}}\right)\widehat{\theta}\left(\frac{\xi_d
}{2^k}-\frac{17}{24}\right)\\
=&~\Phi_k^{+,+}-\Phi_k^{+,-}.
\end{align*}
Similarly,
\begin{align*}
\Phi_k^2=&~\frac{i}{\varepsilon^{2}2^{dk+1}}e^{-i 2^{2n+k}\ee\cdot\left(\xi-\widetilde{\lambda}_n\ee\right)}\widehat{\theta}\left(\frac{\xi_1-\widetilde{\lambda}_n}{\varepsilon2^{k}}\right)
\widehat{\theta}\left(\frac{\xi_2-\widetilde{\lambda}_n}{\varepsilon2^{k}}\right)\widehat{\theta}\left(\frac{\xi_3}{2^{k}}\right)\cdots
\widehat{\theta}\left(\frac{\xi_{d-1}}{2^{k}}\right)\widehat{\theta}\left(\frac{\xi_d
}{2^k}+\frac{17}{24}\right)\\
-&~\frac{i}{\varepsilon^{2}2^{dk+1}}e^{-i 2^{2n+k}\ee\cdot\left(\xi+\widetilde{\lambda}_n\ee\right)}\widehat{\theta}\left(\frac{\xi_1-\widetilde{\lambda}_n}{\varepsilon2^{k}}\right)
\widehat{\theta}\left(\frac{\xi_2-\widetilde{\lambda}_n}{\varepsilon2^{k}}\right)\widehat{\theta}\left(\frac{\xi_3}{2^{k}}\right)\cdots
\widehat{\theta}\left(\frac{\xi_{d-1}}{2^{k}}\right)\widehat{\theta}\left(\frac{\xi_d
}{2^k}-\frac{17}{24}\right)\\
=&~\Phi_k^{-,+}-\Phi_k^{-,-}.
\end{align*}
Recalling that the support condition of $\widehat{\theta}$, we have
\bal\label{f1}
\mathrm{supp} \ \Phi_k^{\pm,\pm}\subset \Big\{\xi:&  \left|\xi_i\pm\frac{17\sqrt{2}}{24}2^{n}\right|\leq \frac{\varepsilon2^k}{100d},\ i=1,2,\nonumber\\
& |\xi_j|\leq \frac{2^k}{100d},\ j=3,\cdots,d-1,\ \left|\xi_d\pm\frac{17}{24}2^k\right|\leq \frac{2^k}{100d}\Big\}.
\end{align}
Without loss of generality, we assume that $\mathrm{supp} \ \widehat{\Phi_k}\subset \mathrm{supp} \ \Phi_k^{+,+}$. Then for all $k\in \mathbb{N}(n)$, we have
\bbal
\mathrm{supp} \ \widehat{\Phi_k}\subset \Big\{\xi:\ \frac{17\sqrt{2}}{24}2^{n}-\frac{\varepsilon2^k}{100d}\leq |\xi_i|\leq \frac{17\sqrt{2}}{24}2^{n}+\frac{\varepsilon2^k}{100d},\ i=1,2,\\
\ |\xi_j|\leq \frac{2^k}{100d},\ j=3,\cdots,d-1,
\ \frac{2}{3}2^k\leq|\xi_d|\leq \frac{3}{4}2^k\Big\},
\end{align*}
which implies that
\bbal
\mathrm{supp} \ \widehat{\Phi_k}\subset  \left\{\xi: \ \frac{33}{24}2^{n}\leq |\xi|\leq \frac{35}{24}2^{n}\right\}.
\end{align*}
Thus, we finish the proof of Lemma \ref{A1}.

\begin{lemma}\label{A2}
Let $H_2$ be defined by \eqref{h2}. Then there holds for $\ell \in \mathbb{N}(n)$
\bal\label{zz}
\dot{\Delta}_\ell H_2=0.
\end{align}
\end{lemma}
{\bf Proof.}\; Obviously, one has
\bbal
H_2&=\frac{1}{2}\sum\limits_{k,j\in \mathbb{N}(n),k\neq j}2^{k+j}\Phi_{k,j}(x)+\frac{1}{2}\sum\limits_{k,j\in \mathbb{N}(n),k\neq j}2^{k+j}\Phi_{k,j}(x)\cos\left(\frac{17}{12}2^{n+1}\ee\cdot x\right)
\end{align*}
where
$$\Phi_{k,j}(x):=\phi\big(2^{k}\e(x-2^{2n+k}\ee)\big)
\phi\big(2^{j}\e(x-2^{2n+j}\ee)\big).$$
Notice that the definition of $\phi$, we deduce that for $j<k$
\bbal
\mathrm{supp} \ \widehat{\Phi_{k,j}}\subset   \left\{\xi\in\R^d: \ \frac{33}{48}2^{k}\leq |\xi|\leq \frac{35}{48}2^{k}\right\},
\end{align*}
which in turn gives that for $j<k$
\bbal
\mathrm{supp} \ \mathcal{F}\left[\Phi_{k,j}\cos\left(\frac{17}{12}2^{n+1}\ee\cdot x\right)\right]\subset   \left\{\xi\in\R^d: \ \frac{33}{24}2^{n+1}\leq |\xi|\leq \frac{35}{24}2^{n+1}\right\}.
\end{align*}
Then, for $j<k$, \eqref{zz} holds. Similarly, \eqref{zz} also holds for $j>k$.

Thus, we finish the proof of Lemma \ref{A2}.
\section*{Acknowledgments}
J. Li is supported by the National Natural Science Foundation of China (11801090 and 12161004) and Jiangxi Provincial Natural Science Foundation (20212BAB211004). Y. Yu is supported by the National Natural Science Foundation of China (12101011) and Natural Science Foundation of Anhui Province (1908085QA05). W. Zhu is supported by the Guangdong Basic and Applied Basic Research Foundation (2021A1515111018).

\section*{Conflict of interest}
The authors declare that they have no conflict of interest.

\end{document}